\documentclass[11pt]{article}
\usepackage{epsf}
\usepackage{latexsym}
\usepackage{amsfonts}

\newtheorem{theorem}{Theorem}[section]
\newtheorem{corollary}[theorem]{Corollary}
\newtheorem{lemma}[theorem]{Lemma}
\newtheorem{definition}[theorem]{Definition}

\def\Esin{\sin_e}
\def\Ecos{\cos_e}
\def\Etan{\tan_e}
\def\Tsin{\sin_t}
\def\Tcos{\cos_t}

\begin{document}
\setlength{\baselineskip}{15pt}

\begin{titlepage}

\title{\bfseries Taxicab Angles and Trigonometry}

\author{Kevin Thompson
	  \thanks{Current address: North Arkansas College,
		1515 Pioneer Drive, Harrison, AR  72601 / kthompson@northark.edu}
	\ \& Tevian Dray
	  \thanks{Send enquiries to: {\tt tevian{\rm @}math.orst.edu}}\\
	Department of Mathematics, Oregon State University,
		Corvallis, OR  97331 \\
}

\date{May 15, 1999 (revised November 5, 1999)}

\thispagestyle{empty}
\renewcommand\thispagestyle[1]{} 
\maketitle

\begin{abstract}
A natural analogue to angles and trigonometry is developed in taxicab
geometry.  This structure is then analyzed to see which, if any, congruent
triangle relations hold.  A nice application involving the use of parallax to
determine the exact (taxicab) distance to an object is also discussed.
\end{abstract}

\end{titlepage}

\section{\textbf{INTRODUCTION}}

Taxicab geometry, as its name might imply, is essentially the study of an
ideal city with all roads running horizontal or vertical.  The roads must be
used to get from point A to point B; thus, the normal Euclidean distance
function in the plane needs to be modified.  The shortest distance from the
origin to the point (1,1) is now 2 rather than $\sqrt{2}$.  So, taxicab
geometry is the study of the geometry consisting of Euclidean points, lines,
and angles in $\mathbb{R}^2$ with the taxicab metric
\[
d((x_{1},y_{1}),(x_{2},y_{2}))=|x_{2}-x_{1}|+|y_{2}-y_{1}|.
\]
A nice discussion of taxicab geometry was given by Krause
\cite{Krause73,Krause75}, and some of its properties have been discussed
elsewhere, including taxicab conic sections \cite{Reynolds,MK,Iny,Laatsch},
and the taxicab isometry group \cite{Schatt}.

In this paper we will explore a slightly modified version of taxicab geometry.
Instead of using Euclidean angles measured in radians, we will mirror the
usual definition of the radian to obtain a taxicab radian (a t-radian).  (A
similar approach was used by Euler \cite{Euler} to discuss the value of $\pi$
for a class of generalized circles which includes the taxicab circle.)  Using
this definition, we will define taxicab trigonometric functions and explore
the structure of the addition formulas from trigonometry.  As applications of
this new type of angle measurement, we will explore the existence of congruent
triangle relations and illustrate how to determine the distance to a nearby
object by performing a parallax measurement.

Henceforth, the label taxicab geometry will be used for this modified taxicab
geometry; a subscript {\it e} will be attached to any Euclidean function or
quantity.

\section{\textbf{TAXICAB ANGLES}}

\begin{figure}[t]
\epsfysize=2.2in
\centerline{\epsffile{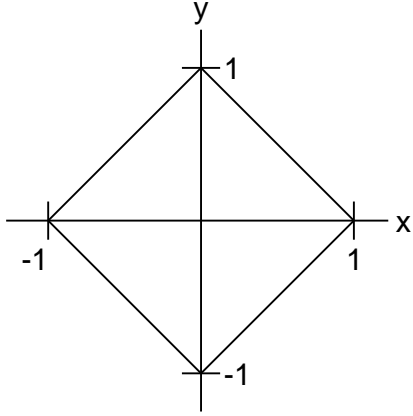}}
\caption{The taxicab unit circle.}
\label{unitcircle}
\end{figure}

There are at least two common ways of defining angle measurement: in terms of
an inner product and in terms of the unit circle.  For Euclidean space, these
definitions agree.  However, the taxicab metric is not an inner product since
the natural norm derived from the metric does not satisfy the parallelogram
law.  Thus, we will define angle measurement on the unit taxicab circle which
is shown in Figure \ref{unitcircle}.

\begin{definition}
A \textbf{t-radian} is an angle whose vertex is the center of a unit (taxicab)
circle and intercepts an arc of (taxicab) length 1.  The taxicab measure of a
taxicab angle $\theta$ is the number of t-radians subtended by the angle on
the unit taxicab circle about the vertex.
\end{definition}

It follows immediately that a taxicab unit circle has 8 t-radians since the
taxicab unit circle has a circumference of 8.  For reference purposes the
Euclidean angles $\pi$/4, $\pi$/2, and $\pi$ in standard position now have
measure 1, 2, and 4, respectively.  The following theorem gives the formula
for determining the taxicab measures of some other Euclidean angles.

\begin{theorem}
An acute Euclidean angle $\phi_e$ in standard position has a taxicab measure
of
\[
\theta
  = 2-\frac{2}{1+\Etan{\phi_e}}
  =\frac{2\Esin{\phi_e}}{\Esin{\phi_e}+\Ecos{\phi_e}}
\]
\label{measure}
\end{theorem}
\textit{Proof:} The taxicab measure $\theta$ of the Euclidean angle $\phi_e$
is equal to the taxicab distance from (1,0) to the intersection of the lines
$y=-x+1$ and $y=x\Etan{\phi_e}$.  The x-coordinate of this intersection is
\[
x_0=\frac{1}{1+\Etan{\phi_e}},
\]
and thus the y-coordinate of P is $y_0=-x_0+1$.  Hence, the taxicab distance
 from (1,0) to P is
\[
\theta=1-x_0+y_0=2-\frac{2}{1+\Etan{\phi_e}}.
\quad\Box
\]

\begin{definition}
The reference angle of an angle $\phi$ is the smallest angle between $\phi$
and the x-axis.
\end{definition}

Theorem \ref{measure} can easily be extended to any acute angle lying entirely
in a quadrant.

\begin{corollary}
If an acute Euclidean angle $\phi_e$ with Euclidean reference angle $\psi_e$
is contained entirely in a quadrant, then the angle has a taxicab measure of
\begin{eqnarray}
\theta
  & = & \frac{2}{1+\Etan{\psi_e}}-\frac{2}{1+\Etan(\phi_e+\psi_e)}
	\nonumber \\
  & = & \frac{2\Esin{\phi_e}}
	{(\Ecos(\phi_e+\psi_e)+\Esin(\phi_e+\psi_e))
	 (\Ecos\psi_e+\Esin\psi_e)} \nonumber
\end{eqnarray}
\label{angleconv}
\end{corollary}

This corollary implies the taxicab measure of a Euclidean angle in
non-standard position is not necessarily equal to the taxicab measure of the
same Euclidean angle in standard position.  Thus, although angles are
translation invariant, they are not rotation invariant.  This is an important
consideration when dealing with any triangles in taxicab geometry.

In Euclidean geometry, a device such as a cross staff or sextant can be used
to measure the angular separation between two objects.  The characteristics of
a similar device to measure taxicab angles would be very strange to
inhabitants of a Euclidean geometry; measuring the taxicab size of the same
Euclidean angle in different directions would usually yield different results.
Thus, to a Euclidean observer the taxicab angle measuring device must
fundamentally change as it is pointed in different directions.  Of course this
is very odd to us since our own angle measuring devices do not appear to
change as we point them in different directions.

The taxicab measure of other Euclidean angles can also be found.  Except for a
few cases, these formulas will be more complicated since angles lying in two
or more quadrants encompass corners of the unit circle.

\begin{figure}[t]
\centerline{\epsffile{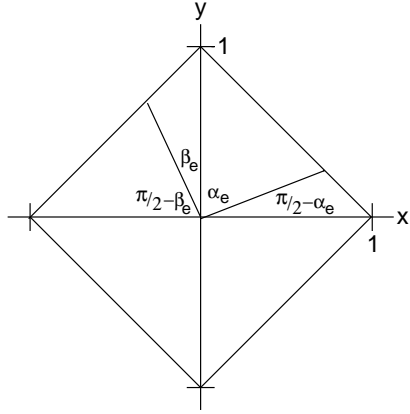}}
\caption{Taxicab right angles are precisely Euclidean right angles.}
\label{right}
\end{figure}

\begin{lemma}
The taxicab measure of any Euclidean right angle is 2 t-radians.
\end{lemma}
\textit{Proof:} Without loss of generality, let $\theta$ be an angle
encompassing the positive y-axis.  As shown in Figure \ref{right}, split
$\theta$ into two Euclidean angles $\alpha_e$ and $\beta_e$ with reference
angles $\pi /2-\alpha_e$ and $\pi /2-\beta_e$, respectively.  Using Theorem
\ref{measure}, we see that
\begin{eqnarray}
\theta
  & = & \frac{2\Esin\alpha_e}
	{\Ecos\alpha_e + \Esin\alpha_e}
	+ \frac{2\Esin\beta_e}{\Ecos\beta_e + \Esin\beta_e} \nonumber \\
  & = & \frac{2\Esin\alpha_e}
	{\Ecos\alpha_e + \Esin\alpha_e}
	+ \frac{2\Ecos\alpha_e}{\Esin\alpha_e + \Ecos\alpha_e} \nonumber \\
  & = & 2 \nonumber
\end{eqnarray}
since $\alpha_e+\beta_e=\pi /2.
\quad\Box$

We now state the taxicab version of the familiar result for the length of an
arc from Euclidean geometry and note that its proof is obvious since all
distances along a taxicab circle are scaled equally as the radius is changed.
This result will be used when we turn to congruent triangle relations and the
concept of parallax.

\begin{theorem}
The length $s$ of the arc intercepted on a (taxicab) circle of radius $r$ by
the central angle with taxicab measure $\theta$ is given by $s=r\theta $.
\label{srtheta}
\end{theorem}

 From the previous theorem we can easily deduce the following result.

\begin{corollary}
Every taxicab circle has 8 t-radians.
\end{corollary}

\section{\textbf{TAXICAB TRIGONOMETRY.}}

We now turn to the definition of the trigonometric functions sine and cosine
in taxicab geometry.  From these definitions familiar formulas for the
tangent, secant, cosecant, and cotangent functions can be defined and results
similar to those below can be obtained.

\begin{definition}
The point of intersection of the terminal side of a taxicab angle $\theta$
in standard position with the taxicab unit circle is the point
$(\Tcos {\theta },\Tsin {\theta })$.
\end{definition}

It is important to note that the taxicab sine and cosine values of a taxicab
angle do not agree with the Euclidean sine and cosine values of the
corresponding Euclidean angle.  For example, the angle 1 t-radian has equal
taxicab sine and cosine values of 0.5.  The range of the cosine and sine
functions remains $[-1,1]$, but the period of these fundamental functions is
now 8.  It also follows immediately (from the distance function) that $|\Tsin
\theta |+|\Tcos \theta |=1$.  In addition, the values of cosine and sine vary
(piecewise) linearly with $\theta$:
\[
\Tcos \theta =\left\{
\begin{array}{l}
1-\frac{1}{2}\theta \;,\;\;\;0\leq \theta <4 \\
\\
-3+\frac{1}{2}\theta \;,\;\;\;4\leq \theta <8
\end{array}
\right.  \;,\;\;\;\;\Tsin \theta =\left\{
\begin{array}{l}
\frac{1}{2}\theta \;,\;\;\;0\leq \theta <2 \\
\\
2-\frac{1}{2}\theta \;,\;\;\;2\leq \theta <6 \\
\\
-4+\frac{1}{2}\theta \;,\;\;\;6\leq \theta <8
\end{array}
\right.
\]

\begin{figure}[t]
\centerline{\epsffile{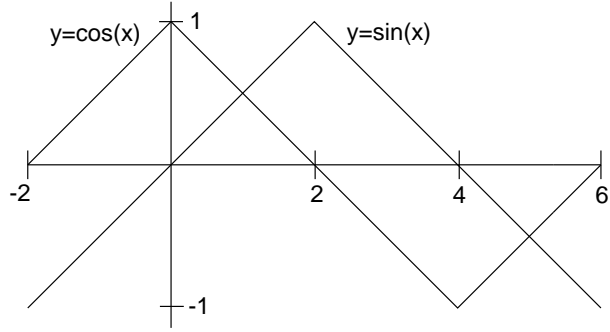}}
\caption{Graphs of the taxicab sine and cosine functions.}
\label{sincos}
\end{figure}

Table \ref{trigrel} gives useful straightforward relations readily derived
 from the graphs of the sine and cosine functions which are shown in Figure
\ref{sincos}.  The structure of the graphs of these functions is similar to
that of the Euclidean graphs of sine and cosine.  Note that the smooth
transition from increasing to decreasing at the extrema has been replaced with
a corner.  This is the same effect seen when comparing Euclidean circles with
taxicab circles.

\begin{table}[b] \centering
\begin{tabular}{|c|c|}
\hline
$\Tsin (-\theta )=-\Tsin \theta $ & $\Tsin (\theta+2)=\Tcos \theta $ \\ \hline
$\Tcos (-\theta )=\Tcos \theta $ & $\Tcos (\theta-2)=\Tsin \theta $ \\ \hline
$\Tsin (\theta-4)=-\Tsin \theta $ &
	$\Tsin (\theta +8k)=\Tsin \theta ,k\in \mathbb{Z}$ \\
\hline
$\Tcos (\theta-4)
  = -\Tcos \theta $ & $\Tcos (\theta +8k)
  = \Tcos \theta ,k\in \mathbb{Z}$
\\ \hline
\end{tabular}
\caption{Basic Taxicab Trigonometric Relations}
\label{trigrel}
\end{table}

As discussed below, and just as in the standard taxicab geometry described in
\cite{Krause75}, SAS congruence for triangles does not hold in modified
taxicab geometry.  Thus, the routine proofs of sum and difference formulas are
not so routine in this geometry.  The first result we will prove is for the
cosine of the sum of two angles.  The formula given for the cosine of the sum
of two angles only takes on two forms; the form used in a given situation
depends on the locations of $\alpha $ and $\beta $.  The notation $\alpha \in
I$ will be used to indicate $\alpha $ is an angle in quadrant I and similarly
for quadrants II, III, and IV.

\begin{theorem}
$\label{cosine}\Tcos (\alpha +\beta )=\pm (-1+|\Tcos \alpha \pm \Tcos \beta |)$
where the signs are chosen to be negative when $\alpha$ and $\beta$ are on
different sides of the \textit{x}-axis and positive otherwise.
\end{theorem}
\textit{Proof:} Without loss of generality, assume $\alpha ,\beta \in [0,8)$,
for if an angle $\theta $ lies outside [0,8), $\exists k\in \mathbb{Z}$ such
that $(\theta +8k)\in [0,8)$ and use of the identity
$\Tcos (\theta +8k)=\Tcos\theta$
will yield the desired result upon use of the following proof.

All of the subcases have a similar structure.  We will prove the subcase
$\alpha\in II, \beta\in III$.  In this situation $6\leq\alpha+\beta\leq10$ and
we take the negative signs on the right-hand side of the equation.  Thus,
\begin{eqnarray}
1-|\Tcos\alpha-\Tcos\beta|
  & = & 1-|1-\frac{1}{2}\alpha-(-3+\frac{1}{2}\beta)| \nonumber \\
  & = & 1-|4-\frac{1}{2}(\alpha+\beta)| \nonumber
\end{eqnarray}
$$
\hspace{1.5in} \quad\qquad
  =\left\{
\begin{array}{l}
-3+\frac{1}{2}(\alpha+\beta) \;\; ,
	\quad 6\leq\alpha+\beta<8 \\ \\
	5-\frac{1}{2}(\alpha+\beta) \; , \qquad 8\leq\alpha+\beta\leq 10
\end{array}
\right.
$$

$$
\hspace{2in} \!= \quad\Tcos(\alpha+\beta) \hspace{1.55in} \Box
$$

\begin{corollary}
$\Tcos (2\alpha )=-1+2|\Tcos \alpha |$.
\end{corollary}

The curious case structure in Theorem \ref{cosine} is due to the odd
combinations of quadrants that determine which sign to choose.  The reason for
the sign change when $\alpha $ and $\beta $ are on different sides of the
\textit{x}-axis lies in the fact that a corner of the cosine function is being
crossed (i.e.  different pieces of the cosine function are being used) to
obtain the values of the cosine of $\alpha $ and $\beta $.  Table
\ref{cosregion} summarizes which form of $\Tcos(\alpha+\beta)$ should be used
when.

\begin{table}[h] \centering
\begin{tabular}{|c|cc|}
\hline
& $\quad \alpha \quad $ & $\quad \beta \quad $ \\ \hline
& same & quadrant \\
& I & II\\
\raisebox{3ex}[0pt]
{$\Tcos(\alpha+\beta) = -1+|\Tcos\alpha+\Tcos\beta|$}& III & IV\\ \hline
& I & III\\
& I & IV\\
& II & III\\
\raisebox{4.5ex}[0pt]
{$\Tcos(\alpha+\beta) = 1-|\Tcos\alpha-\Tcos\beta|$} & II & IV\\ \hline
\end{tabular}
\caption{Forms of $\Tcos(\alpha+\beta)$ and Regions of Validity}
\label{cosregion}
\end{table}

\newpage

We can use Theorem \ref{cosine} and the relations in Table \ref{trigrel} to
establish a pair of corollaries.

\begin{corollary}
$\Tsin (\alpha +\beta )=\pm (-1+|\Tsin \alpha \pm \Tcos \beta |)$ where the
signs are chosen according to Table \ref{sinregion}.
\label{sine}
\end{corollary}
\textit{Proof:} First, note $\Tsin\theta=\Tcos(\theta-2)$.  As with the cosine
addition formula, all cases are proved similarly.  We will assume $\alpha \in
I$ and $\beta \in IV$.  We have $\alpha-2$ and $\beta$ in the same quadrant,
and thus
\begin{eqnarray}
\Tsin(\alpha+\beta) & = & \Tcos((\alpha+\beta)-2) \nonumber \\
                   & = & \Tcos((\alpha-2)+\beta) \nonumber \\
                   & = & -1 + |\Tcos(\alpha-2)+\Tcos \beta| \nonumber \\
	           & = & -1 + |\Tsin\alpha+\Tcos \beta|
\quad \quad \Box \nonumber
\end{eqnarray}

\begin{corollary}
$\Tsin (2\alpha )=-1+2|\Tcos(\alpha -1)|$
\end{corollary}
\textit{Proof:}
$\Tsin(2\alpha)=\Tcos(2\alpha-2)=\Tcos(2(\alpha-1))=-1+2|\Tcos(\alpha-1)|
\quad \quad \qquad \Box$

\begin{table}[t] \centering
\begin{tabular}{|c|cc|}
\hline
& $\quad \alpha \quad $ & $\quad \beta \quad $ \\ \hline
& I & III\\
& I & IV\\
& II & II\\
\raisebox{4.5ex}[0pt]
{$\Tsin(\alpha+\beta)=-1+|\Tsin\alpha+\Tcos\beta|$}& IV & IV\\ \hline
& I & I\\
& I & II\\
& II & III\\
& II & IV\\
& III & III\\
\raisebox{7ex}[0pt]
{$\Tsin(\alpha+\beta)=1-|\Tsin\alpha-\Tcos\beta|$}& III & IV\\ \hline
\end{tabular}
\caption{Forms of $\Tsin(\alpha+\beta)$ and Regions of Validity}
\label{sinregion}
\end{table}

\section{\textbf{CONGRUENT TRIANGLES}}

\begin{figure}[b]
\centerline{\epsffile{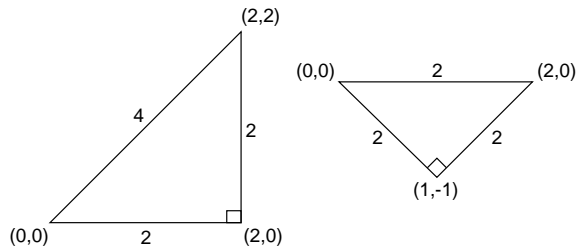}}
\caption{Triangles satisfying ASASA that are not congruent.}
\label{triangles}
\end{figure}

In Euclidean geometry we have many familiar conditions that ensure two
triangles are congruent.  Among them are SAS, ASA, and AAS.  In modified
taxicab geometry the only condition that ensures two triangles are congruent
is SASAS.  One example eliminates almost all of the other conditions.

Consider the two triangles shown in Figure \ref{triangles}.  The triangle
formed by the points (0,0), (2,0), and (2,2) has sides of lengths 2, 2, and 4
and angles of measure 1, 1, and 2 t-radians.  The triangle formed by the
points (0,0), (2,0), and (1,-1) has sides of length 2 and angles of measure 1,
1, and 2.  These two triangles satisfy the ASASA condition but are not
congruent.  This also eliminates the ASA, SAS, and AAS conditions as well the
possibility for a SSA or AAA condition.

The triangle formed by the points (0,0), (0.5,1.5), and (1.5, 0.5) has sides
of length 2 and angles 1, 1.5, and 1.5 t-radians.  Thus, it satisfies the SSS
condition with the second triangle in the previous example.  However, the
angles of these triangles and not congruent.  Hence, the SSS and SSSA
conditions fail.

The last remaining condition, SASAS, actually does hold.  Its proof relies on
the fact that even in this geometry the sum of the angles of a triangle is a
constant 4 t-radians, which in turn relies on the fact that, given parallel
lines and a transversal, alternate interior angles are congruent.  We begin by
noting that opposite angles are congruent.  This leads immediately to the
following result.

\begin{figure}[t]
\epsfysize=1.1in
\centerline{\epsffile{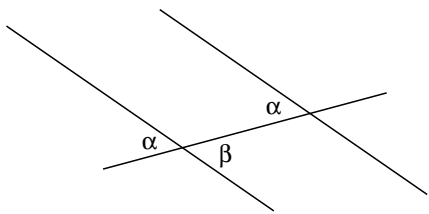}}
\caption{Alternate interior angles formed by parallel lines and a transversal
are congruent.}
\label{altint}
\end{figure}

\begin{lemma}
Given two parallel lines and a transversal, the alternate interior angles are
congruent.
\end{lemma}
\textit{Proof:} Using Figure \ref{altint} translate $\alpha$ along the
transversal to become an angle opposite $\beta$.  By the note above, $\alpha$
and $\beta$ are congruent.
$\qquad \Box$

\begin{figure}[b]
\epsfysize=2in
\centerline{\epsffile{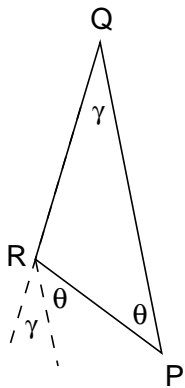}}
\caption{The sum of the angles of a taxicab triangle is always 4 t-radians.}
\label{sumangles}
\end{figure}

\begin{theorem}
The sum of the angles of a triangle in modified taxicab geometry is 4
t-radians.
\end{theorem}
\textit{Proof:} Given the triangle in Figure \ref{sumangles}, we can translate
the angle $\gamma$ from Q to R and by the congruence of alternate
interior angles conclude the sum of the angles of the triangle is 4 t-radians.
$\quad \Box$

Therefore, given two triangles having all three sides and any two angles
congruent, the triangles must be congruent.  However, as we have seen, this is
the only congruent triangle relation in taxicab geometry.

\begin{figure}[t]
\epsfysize=2.8in
\centerline{\epsffile{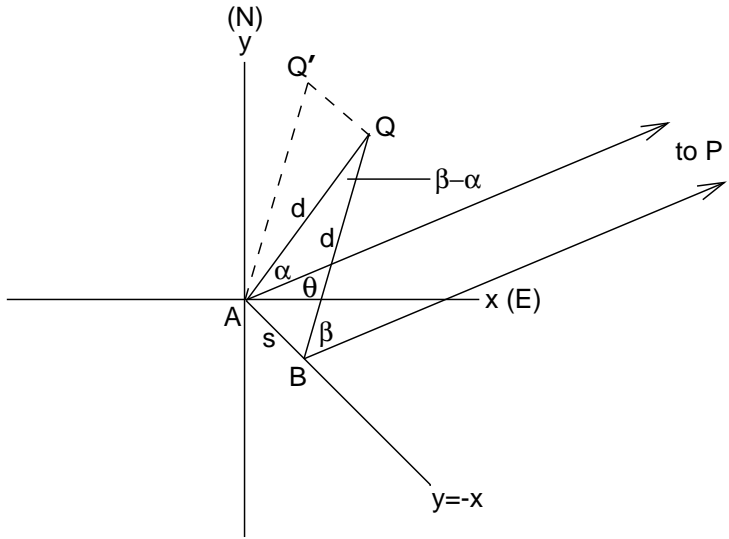}}
\caption{A parallax diagram in taxicab geometry.}
\label{taxipara}
\end{figure}

\section{\textbf{PARALLAX}}

Parallax, the apparent shift of an object due to the motion of the observer,
is a commonly used method for estimating the distance to a nearby object.  The
method of stellar parallax was used extensively to find the distances to
nearby stars in the $19^{th}$ and early $20^{th}$ centuries.  We now wish to
explore the method and results of parallax in taxicab geometry and examine how
these differ from the Euclidean method and results.  We will discover that the
taxicab method yields the same formula commonly used in the Euclidean case
with the exception that the taxicab formula is exact.

Suppose that as a citizen of Modified Taxicab-land you wish to find the
distance to a nearby object $Q$ in the first quadrant, and that there is also
a distant reference object $P$ ``at infinity'' essentially in the same
direction as $Q$ with reference angle $\theta$ (Figure \ref{taxipara}).  The
distant reference object should be far enough away so that it appears
stationary when you move small distances.  We may assume without loss of
generality that the object $Q$ does not lie on either axis, for if it did, we
could move a small distance to get the object in the interior of the first
quadrant.

Initially standing at $A$, measure the angle $\alpha$ between $Q$ and $P$
using the taxicab equivalent of a cross staff or a sextant.  Now, for reasons
to be apparent later, you should move a small distance (relative to the
distance to the object) in such a way that the distance to the object does not
change.  This can be accomplished by moving in either of two directions, and,
provided you move only a small distance, the object remains in the interior of
the first quadrant.  Furthermore, exactly one of these directions results in
the angle between $Q$ and $P$ being increased, so that the situation depicted
in Figure \ref{taxipara} is generic.

You have therefore moved from $A$ to $B$ in one of the following directions:
NW, NE, SW, or SE.  Now measure the new angle $\beta $ between the two
objects.  With this information we can now find the taxicab distance to the
object $Q$.  Construct the point $Q^{\prime }$ such that $\overline{QQ^{\prime
}}$ is parallel to $\overline{AB}$ and $\ell (\overline{QQ^{\prime }})=\ell
(\overline{AB})=s$.  The angle $\angle PAQ^{\prime }$ has measure $\beta $
since it is merely a translation of $\angle QBP$.  Thus, $\angle QAQ^{\prime
}$ has measure $\beta -\alpha $.  Now, the lengths of $\overline{AQ}$ and
$\overline{BQ}$ are equal since the direction of movement from $A$ to $B$ was
shrewdly chosen so that the distance to the object remained constant.  Since
translations do not affect lengths, this implies $AQ$ and $AQ^{\prime }$ have
equal lengths.  Hence, the points $Q$ and $Q^{\prime }$ lie on a taxicab
circle of radius \textit{d} centered at $A$.  Using the formula for the length
of a taxicab arc in Theorem \ref{srtheta},
\begin{equation}
d=\frac{s}{\beta -\alpha }  \label{taxicab}
\end{equation}
where \textit{s} and \textit{d} are taxicab distances and $\beta -\alpha $ is
a taxicab angle.  This formula is identical to the Euclidean distance
\textit{estimation} formula with \textit{d} and \textit{s} Euclidean distances
and $\beta -\alpha $ an Euclidean angle.  However, as we shall now see, the
commonly used Euclidean version is truly an approximation and not an exact
result.  This realization is necessary to logically link the commonly used
Euclidean formula and the taxicab formula.

\begin{figure}[t]
\epsfysize=2.3in
\centerline{\epsffile{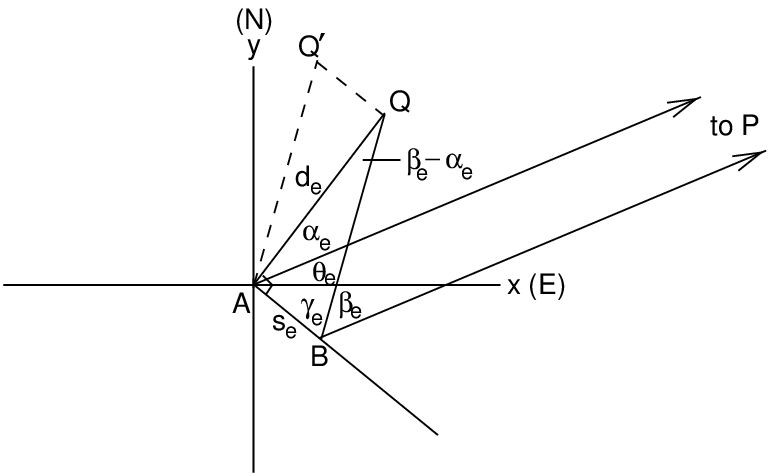}}
\caption{A parallax diagram in Euclidean geometry.}
\label{euclidpara}
\end{figure}

Using Figure \ref{taxipara} but with all distances and angles now Euclidean,
we isolate $\triangle QAB$.  Using the law of sines and the fact that $\gamma
=3\pi /4-(\beta_e +\theta_e )$, we have
\begin{equation}
d_e=\frac{s_e(\Ecos(\beta_e +\theta_e )
  + \Esin(\beta_e+\theta_e))}{ \sqrt{2}\Esin(\beta_e -\alpha_e )}
  \label{euclidean}
\end{equation}
This formula can be simplified by moving from $A$ to $B$ in a direction
perpendicular to the line of sight to Q from A rather than in one of the four
prescribed directions above.  In this case $m(\angle QBA)=\pi /2-(\beta_e
-\alpha_e )$ (Figure \ref{euclidpara}).  Thus, the law of sines gives the
Euclidean parallax formula
\[
d_e=\frac{s_e}{\Etan(\beta_e -\alpha_e )}
\]
If we now apply the approximation $\Etan(\beta_e -\alpha_e )\approx
(\beta_e -\alpha_e )$ for small angles we obtain the commonly used Euclidean
parallax formula
\[
d_e\approx \frac{s_e}{\beta_e -\alpha_e}.
\]
which is not exact.

It is interesting to note the quite different movement requirements in the two
geometries needed to obtain the best possible approximations of the distance
to the object.  This difference lies in the methods of keeping the distance to
the object as constant as possible.  In the Euclidean case, moving small
distances on the line tangent to the circle of radius \textit{d} centered at
the object (i.e.  perpendicular to the radius of this circle) essentially
leaves the distance to the object unchanged.  In taxicab geometry, moving in
one direction along either $y=x$ or $y=-x$ keeps the distance to the object
exactly unchanged.

We are now in a position to justify the link between results (\ref{taxicab})
and (\ref{euclidean}).  Since the line segment $\overline{QQ^{\prime }}$ of
taxicab length \textit{s} lies on a taxicab circle, $s=\sqrt{2}s_e$.  The
distance \textit{d} to the object is given by $d=d_e(\Ecos(\alpha_e
+\theta_e )+\Esin(\alpha_e +\theta_e ))$ since the Euclidean angle between
the line of sight $\overline{AQ}$ and the x-axis is $(\alpha_e +\theta_e )$.
Using Corollary \ref{angleconv} with $\phi =(\beta_e -\alpha_e )$ and $\psi
=(\alpha_e +\theta_e )$, the taxicab measure of $\beta -\alpha $ is given by
\[
\beta -\alpha
  = \frac{2\Esin(\beta_e -\alpha_e )}{(\Ecos(\beta_e +\theta_e)
    + \Esin(\beta_e +\theta_e ))(\Ecos(\alpha_e +\theta_e)
    + \Esin(\alpha_e+\theta_e ))}.
\]
Using these substitutions, formula (\ref{taxicab}) becomes formula
(\ref{euclidean}).

\section{\textbf{CONCLUSION}}

With this natural definition of angles in taxicab geometry, some of the same
difficulties arise as with Euclidean angles in taxicab geometry.  Congruent
triangles are few and far between.  Only with the strictest requirements,
namely that all three sides and two angles are congruent, are we able to
conclude that two triangles must be congruent.

In addition to creating a natural definition of angles and trigonometric
functions, we have also unwittingly created an environment in which a parallax
method can be used to determine the exact distance to a nearby object rather
than just an approximation.  This is not too surprising a result since there
exist directions in which one can travel without the distance to an object
changing.  With this result, and taxicab and Euclidean angle measuring
instruments, the exact Euclidean distance to the object can now be found (up
to measurement error of course).  The trick is to build your own taxicab cross
staff or sextant.

\section*{\textbf{ACKNOWLEDGMENTS}}

This work is based on a paper submitted to Oregon State University in early 1997 by Kevin Thompson in partial fulfillment of the requirements for the M.S.\ degree in Mathematics.

\end{document}